# About one method of constructing Hermite trigonometric splines

Denysiuk V.P. Dr of phys.-math. sciences, Professor, Hryshko O.N.
Kiev, Ukraine
National Aviation University
kvomden@nau.edu.ua

## Annotation

The method of constructing trigonometric Hermite splines, which interpolate the values of some periodic function and its derivatives in the nodes of a uniform grid, is considered. The proposed method is based on the periodicity properties of trigonometric functions and is reduced to solving only systems of linear algebraic equations of the second order; solutions of these systems can be obtained in advance. When implementing this method, it is necessary to calculate the coefficients of interpolation trigonometric polynomials that interpolate the values of the function itself and the values of its derivatives at the nodes of the uniform grid; known fast Fourier transform algorithms can be used for this purpose. Examples of construction of trigonometric Hermite splines of the first and second orders are given. The proposed method can be recommended for practical use.

## Keywords:

Hermite interpolation, trigonometric, Hermite interpolation, trigonometric splines, trigonometric Hermite splines.

## Introduction

In many problems of science and technology, the problem of interpolation of a function in such a statement often arises. Let $C^r$, ($r = 1, 2, \dots$ ) - a class of functions having continuous derivatives up to $r$ - th order inclusively. Let also on the segment $[a, b]$ some grid is set $a = x_0 < x_1 < \dots < x_N = b$ and let the values of the function itself be known in the grid nodes $f(x_i) = f_i$ ($i = 0, 1, \dots, N$), and the meaning of derivatives $f^{(k)}(x_i) = f_i^{(k)}$, ($k = 1, 2$) of this function. You need to build a function $\mathrm{H}(x)$ of class $C^r$, ($r - 1, 2, \dots$ ), depending on $(N+1)(k+1)$ parameters that satisfies the conditions

$$\mathrm{H}^{(k)}(x_i) = f_i^{(k)}, \qquad (i = 0, 1, \dots, N; k = 0, 1, \dots, r). \tag{1}$$

In the case where the function $\mathrm{H}(x) \in \mathrm{P}$, where $\mathrm{P}$ is a class of algebraic polynomials, the polynomial that satisfies these conditions is called the Hermit interpolation polynomial.

In many cases, a uniform grid is given on the interval $[0, 2\pi)$, and the function $\mathrm{H}(x)$ s sought in the class of trigonometric polynomials $\mathrm{T}$; the polynomial $T(x)$ that is the solution of this problem is called the Hermitian trigonometric polynomial. The problem of constructing a trigonometric Hermite polynomial that satisfies conditions (1) was considered in [Den].

A certain disadvantage of the Hermit trigonometric polynomial is its analyticity. Therefore, in our opinion, the problem of constructing Hermit trigonometric splines belonging to classes $C^r$, ($r = 1, 2$); In the future, such polynomials will be called Hermite trigonometric splines.

## Analysis of researches and publications

There are different approaches to constructing Hermite interpolation splines. For example, one of them is based on a modification of Hermite polynomial interpolation methods. Note that the problem of constructing a Hermite polynomial polynomial on a segment in general is considered in detail, for example, in [2]; the construction of Hermite polynomial splines is considered in [3,4].

The construction of Hermite trigonometric polynomials and the problems associated with this approach are considered, for example, in [5]; there is also a sufficient bibliography on this issue. This approach and its development are also considered in [6], [7].

Another approach to the construction of Hermite trigonometric polynomials is based on the application of fundamental trigonometric functions. This approach is considered, for example, in [8]; it is expedient to emphasize that in this work the fundamental trigonometric functions are constructed, in our opinion, in the most natural way. Other variants of this approach, which differ in the methods of constructing fundamental functions, are considered in [9], [10], [11]. The disadvantage of this approach is the difficulty of constructing fundamental trigonometric functions of higher orders; the consequence of this is that the constructed Hermite polynomials interpolate the function itself and its derivative only of the first order.

Finally, another approach to the construction of trigonometric Hermite polynomials of arbitrary order is given in [12]. This approach is based on the periodicity properties of trigonometric functions and is reduced to solving only two systems of linear algebraic equations $k+1$ -th order; solutions of these systems can be obtained in advance. When implementing this approach, it is necessary to calculate the coefficients of interpolation trigonometric polynomials that interpolate the values of the function itself and the values of its derivatives in the nodes of a uniform grid; known fast Fourier transform algorithms can be used for this purpose. There are also examples of constructing trigonometric Hermite polynomials for $p = 1,2$ .

As we have already said, a certain disadvantage of the Hermite trigonometric polynomial is its analyticity. Therefore, in our opinion, the problem of constructing trigonometric Hermite's functions belonging to classes is interesting. $C^r$, ($r = 1,2$); in the future such functions will be called Hermite trigonometric splines.

## The purpose of the work.

The aim of this work is to develop a constructive method for constructing Hermite trigonometric splines that interpolate a function and its derivatives of order $p$ , ($p = 1,2$), in nodes of uniform grids on an interval $[0, 2\pi)$, and belong to spaces $C^r$, (r=1,2).

## The main part.

During interpolating functions, it is often the case that in the nodes of a uniform grid $\Delta_N^l$ , given on the segment $[0, 2\pi]$ , not only the values of the function $f(t)$ are known, but also the values of the derivatives of this function. In the case when the interpolation is carried out by algebraic polynomials, Hermite's polynomials are used, which interpolate in the nodes of the grid $\Delta_N^l$ both the value of the function itself and the value of its derivatives. If the interpolation of the function is supposed to be carried out by trigonometric polynomials, then the question arises about the construction of trigonometric polynomials, which are analogs of algebraic Hermite polynomials; such trigonometric polynomials will hereinafter be called Hermite trigonometric polynomials.

As we have already said, a certain disadvantage of the Hermite trigonometric polynomial is its analyticity. Therefore, in our opinion, the problem of constructing trigonometric Hermite functions belonging to classes $C^r$, ($r = 1,2$) is interesting; in the future such functions will be called Hermite trigonometric splines. Consider the construction of trigonometric Hermite splines in more detail, starting with the case $r = 1$ .

Let on $[0, 2\pi)$ be given a grid $\Delta_N^{(I)} = \left\{ t_j^{(I)} \right\}_{j=1}^N$ , ($I = 0,1$; $N = 2n+1, n = 1,2,...$; $t_j^{(0)} = \frac{2\pi}{N}(j-1)$ ; $t_j^{(1)} = \frac{\pi}{N}(2j-1)$ ), and let the nodes of this grid be given the value of the function $f(t)$ and its first-order derivatives. Denote the value of the function and its derivatives in the nodes of the grid $\Delta_N^{(I)}$ by $f_j^{(I,0)}, f_j^{(I,1)}$, ($j = 1,2,...,N$). With respect to these derivatives, we will assume that these derivatives are centered on the grid $\Delta_N^{(I)}$, the condition is met

$$\sum_{j=1}^N f_j^{(I,1)} = 0.$$

Let's put $I = 0$ and we look for a trigonometric Hermite spline that interpolates a function and its first-order derivatives at grid nodes $\Delta_N^{(0)}, (N = 2n+1, n = 1,2,...)$, like

$$Se_{1,n}^{(0)}(t) = \frac{a_0^{(0)}}{2} + \sum_{k=1}^{n}\left[a0_k^{(0)}C0^{(0)}(k,t) + b0_k^{(0)}S0^{(0)}(k,t) + a1_k^{(0)}C1^{(0)}(k,t) + b1_k^{(0)}S1^{(0)}(k,t)\right], \qquad (1)$$

where

$$C0^{(0)}(k,t) = \frac{\cos(kt)}{k^3} + \sum_{m=1}^{\infty}\frac{\cos((mN+k)t)}{(mN+k)^3}; \quad S0^{(0)}(k,t) = \frac{\sin(kt)}{k^3} + \sum_{m=1}^{\infty}\frac{\sin((mN+k)t)}{(mN+k)^3};$$

$$C1^{(0)}(k,t) = \sum_{m=1}^{\infty}\frac{\cos((mN-k)t)}{(mN-k)^3}; \qquad S1^{(0)}(k,t) = \sum_{m=1}^{\infty}\frac{\sin((mN-k)t)}{(mN-k)^3},$$

and $a_0^{(0)}$, $a0_j^{(0)}$, $b0_j^{(0)}$, $a1_j^{(0)}$, $b1_j^{(0)}$ ( $j = 1,2,...,n$ ) - are uncertain coefficients.

It is clear that the series that define the function $C0^{(0)}(k,t)$, $S0^{(0)}(k,t)$, $C1^{(0)}$, $S1^{(0)}(k,t)$ on the basis of Weierstrass are uniformly convergent series, therefore, they provide a continuous function. Further, these series can be differentiated by members, and the first derivative of these functions will also be continuous. Thus, it is clear that all functions $C0^{(0)}(k,t)$, $S0^{(0)}(k,t)$, $C1^{(0)}$, $S1^{(0)}(k,t)$, and hence the spline itself $Se_{1,n}^{(0)}(t)$ belong to the class $C_{[0,2\pi)}^1$.

To find the coefficients $a0_j^{(0)}$, $b0_j^{(0)}$, $a1_j^{(0)}$, $b1_j^{(0)}$ let's do this. Calculate the value of the spline $Se_{1,n}^{(0)}(t)$ in grid nodes $\Delta_N^{(0)}$ and equate them to the values of the trigonometric polynomial $T0_n^{(0)}(t)$, interpolating the value of the function $f_j^{(0,0)}$ in nodes of the same grid. Next, differentiating the spline $Se_{1,n}^{(0)}(t)$ and equating the value of the derivative at the nodes of the grid $\Delta_N^{(0)}$ the value of the interpolation trigonometric polynomial $T1_n^{(0)}(t)$, interpolating the value of the same derivative function $f_j^{(0,1)}$ in grid nodes $\Delta_N^{(0)}$, we receive $2n$ systems of linear algebraic equations 2 systems of linear algebraic equations $Se_{1,n}^{(0)}(t)$.

The proposed approach to the construction of trigonometric Hermite splines is based on the fact that the values of trigonometric functions $\cos(mN \pm k)t$ and $\sin(mN \pm k)t$, ( $m = 1,2,...$ ) in grid nodes $\Delta_N^{(0)}$ depend only on the node numbers of this grid $j$ and indices $k$; that is, equality takes place $\cos(mN+k)t_j^{(0)} = \left|\cos kt_j^{(0)}\right|$, $\sin(mN \pm k)t_j^{(0)} = \left|\sin kt_j^{(0)}\right|$. Using this fact can greatly simplify the task of finding the coefficients $a_j^{(0)}$, $b_j^{(0)}$, spline $Se_{1,n}^{(0)}(t)$, reducing it to the problem of calculating the coefficients of interpolation trigonometric polynomials and the subsequent solution $N-1$ systems of algebraic equations of dimension 2.

To illustrate, we construct a trigonometric Hermite spline that interpolates the function $f(t)$ and its derivative 1 order in grid nodes $\Delta_N^{(0)}$, which we will look for in the form (1); at the same time for reduction we will lower a grid index 0.

Calculate the value of the spline $Se_{1,n}(t)$ in grid nodes $\Delta_N$. Taking into account that $t_i = \frac{2\pi}{N}(j-1)$ and equality

$$\cos(N-k)t_j = \cos(N-k)\frac{2\pi}{N}(j-1) = \cos kt_j,$$

$$\sin(N-k)t_j = \sin(N-k)\frac{2\pi}{N}(j-1) = -\sin kt_j,$$

we have

$$Se_{1,n}(t) = \frac{a_0}{2} + \sum_{k=1}^{n}\left[(gc0(k)a0_k + gc1(k)a1_k)\cos kt_j + (gs0(k)b0_k - gs1(k)b1_k)\sin kt_j\right], \quad (2)$$

where

$$gc0(k) = \frac{1}{k^3} + \sum_{m=1}^{\infty}\frac{1}{(mN+k)^3}; \quad gs0(k) = \frac{1}{k^3} + \sum_{m=1}^{\infty}\frac{1}{(mN+k)^3};$$

$$gc1(k) = \sum_{m=1}^{\infty} \frac{1}{(mN-k)^3}; \qquad gs1(k) = -\sum_{m=1}^{\infty} \frac{1}{(mN-k)^3},$$

Equating the obtained values of the spline $Se_{1,n}(t)$ in grid nodes $\Delta_N$ the value of the trigonometric interpolation polynomial $T0_n(t)$, interpolating the value of the function $f_j^{(0,0)}$ in nodes of the same grid, we obtain relative to the coefficients $a_k^{(0)}$, $a_{N-k}^{(0)}$, $b_k^{(0)}$, $b_{N-k}^{(0)}$ such systems of equations of the second order

$$a0_0 = A0_0;$$

$$\begin{cases} gc0(k)a0_k + gc1a1_k = A0_k \\ gs0(k)b_k - gs1(k)b_{N-k} = B0_k; \end{cases} \qquad k=1,2,...,n. \tag{3}$$

where $A0_0$, $A0_k$ and $B0_k$ - coefficients of the trigonometric polynomial $T0_n(t)$.

Differentiating spline $Se_{1,n}(t)$, we have

$$Se'_{1,n}(t) = \sum_{k=1}^{n} \Big[ a0_k C01(k,t) + b0_k S01(k,t) + a1_k C11(k,t) + b1_k S11(k,t) \Big],$$

where

$$C01(k,t) = -\frac{\sin(kt)}{k^2} - \sum_{m=1}^{\infty} \frac{\sin((mN+k)t)}{(mN+k)^2}; \quad S01(k,t) = \frac{\cos(kt)}{k^2} + \sum_{m=1}^{\infty} \frac{\cos((mN+k)t)}{(mN+k)^2};$$

$$C11(k,t) = -\sum_{m=1}^{\infty} \frac{\sin((mN-k)t)}{(mN-k)^2}; \qquad S11(k,t) = -\sum_{m=1}^{\infty} \frac{\cos((mN-k)t)}{(mN-k)^2}.$$

Calculating the value of the spline $Se'_{1,n}(t)$ in grid nodes $\Delta_N$, we receive

$$Se'_{1,n}(t_j) = \sum_{k=1}^{n} \Big[ (gs01(k)b0_k + gs11(k)b1_k)\cos(kt_j) + (gc01(k)a0_k + gs11(k)a1_k)\sin(kt_j) \Big].$$

where

$$gc01(k) = -\left[ \frac{1}{k^2} + \sum_{m=1}^{\infty} \frac{1}{(mN+k)^2} \right]; \qquad gs0(k) = \frac{1}{k^2} + \sum_{m=1}^{\infty} \frac{1}{(mN+k)^2};$$

$$gc11(k) = \sum_{m=1}^{\infty} \frac{1}{(mN-k)^2}; \qquad gs1(k) = -\sum_{m=1}^{\infty} \frac{1}{(mN-k)^2},$$

Marking through $A1_k$, $B1_k$ coefficients of a trigonometric polynomial that interpolates the value of the derivative $f'(t)$ in grid nodes $\Delta_N$, and given that $A1_0 = 0$, we obtain systems of equations

$$\begin{cases} gs01(k)b0_k + gs11(k)b1_k = A1_k; \\ gc01(k)a0_k + gs11(k)a1_k = B1_k. \end{cases} \tag{4}$$

Systems of equations (3), (4) can be submitted as

$$a0_0 = A0_0;$$

$$\begin{cases} gc0(k)a0_k + gc1a1_{N-k} = A0_k; \\ gc01(k)a0_k + gs11(k)a1_k = B1_k. \end{cases} ; \qquad \begin{cases} gs0(k)b0_k - gs1(k)b1_k = B0_k \\ gs01(k)b0_k + gs11(k)b1_k = A1_k. \end{cases} ;$$

$$k = 1,2,...,n.$$

Solving these system, we have

$$a0(k) = \frac{1}{Qc(k)}\left( \frac{gc11(k)}{gc1(k)} A0(k) - B1(k) \right); \qquad b0(k) = \frac{1}{Qs(k)}\left( \frac{-gs11(k)}{gs1(k)} B0(k) + A1(k) \right);$$

$$a1(k) = \frac{1}{Qc(k)}\left( \frac{-gc01(k)}{gc0(k)} A0(k) + B1(k) \right); \qquad b1(k) = \frac{1}{Qs(k)}\left( \frac{-gs01(k)}{gs0(k)} B0(k) + A1(k) \right),$$

$$k = 1,2,...,n.$$

where

$$Qc(k) = \frac{gc11(k)}{gc1(k)} - \frac{gc01(k)}{gc0(k)} \; ; \qquad Qs(k) = -\left( \frac{gs11(k)}{gs1(k)} - \frac{gs01(k)}{gs0(k)} \right).$$

Substituting the found coefficients in (1), we obtain a trigonometric Hermite spline $Se_{1,n}(t)$, which is sewn in knots of a grid $\Delta_N^{(0)}$ and interpolates the function $f(t)$ and its derivative at the nodes of the same grid..

Thinking similarly, it is easy to obtain a trigonometric Hermite polynomial that is crosslinked at the nodes of the grid. $\Delta_N^{(1)}$ and interpolates the function $f(t)$ and its derivative at the nodes of the grid $\Delta_N^{1}$; such a spline has the form

$$Se_{2n}^{(1)}(t) = \frac{a_0^{(1)}}{2} + \sum_{k=1}^{n} \left[ a0_k^{(1)} C0^{(0)}(k,t) + b0_k^{(1)} S0^{(0)}(k,t) + a1_k^{(1)} C1^{(0)}(k,t) + b1_k^{(1)} S1^{(0)}(k,t) \right], \qquad (6)$$

where the coefficients $a0_j^{(1)}, b0_j^{(1)}, a1_j^{(1)}, b1_j^{(1)}$ are determined by formulas

$$a0_k = \frac{1}{Qc0(k)} \left( \frac{gc11(k)}{gc1(k)} A0_k - B1_k \right); \qquad b0_k = \frac{1}{Qc0(k)} \left( -\frac{gs11(k)}{gs1(k)} B0_k + A1_k \right);$$

$$a1_k = \frac{1}{Qc0(k)} \left( -\frac{gc01(k)}{gc0(k)} A0_k + B1_k \right); \qquad b1_k = \frac{1}{Qs0(k)} \left( -\frac{gs01(k)}{gs0(k)} B0_k + A1_k \right);$$

$$gc0^{(1)}(k) = gs0^{(1)}(k) = \frac{1}{k^3} + \sum_{m=1}^{\infty} \frac{(-1)^m}{(mN+k)^3} \; ; \quad -gc1^{(1)}(k) = gs1^{(1)}(k) = \sum_{m=1}^{\infty} \frac{(-1)^m}{(mN-k)^3} \; ,$$

$$-gc01^{(1)}(k) = gs01^{(1)}(k) = \frac{1}{k^2} + \sum_{m=1}^{\infty} \frac{(-1)^m}{(mN+k)^2} \; ; \quad gc11^{(1)}(k) = gs11^{(1)}(k) = \sum_{m=1}^{\infty} \frac{(-1)^{m+1}}{(mN+k)^2} \; ;$$

$$Qc^{(1)}(k) = Qs^{(1)}(k) = \frac{gc11^{(1)}(k)}{gc1^{(1)}(k)} - \frac{gc01^{(1)}(k)}{gc0^{(1)}(k)} \; ;$$

and coefficients $A0_0^{(1)}, A0_k^{(1)} \; B0_k^{(1)}, \; A1_k^{(1)}, \; B1_k^{(1)}$ are the coefficients of interpolation trigonometric polynomials $T0_{2n}^{(1)}$ and $T1_{2n}^{(1)}$ interpolating the function and its derivative in the grid nodes, respectively $\Delta_N^{(1)}$.

Let us now consider the construction of a trigonometric Hermite spline that interpolates a function and its first and second derivatives at grid nodes. $\Delta_N^{(I)}$. This spline has the form

$$Se_{2,n}^{(I)}(t) = \frac{a_0^{(I)}}{2} + \sum_{k=1}^{n} \left\{ a0_k^{(I)} C0^{(I)}(k,t) + a1_k^{(I)} C1^{(I)}(k,t) + a2_k^{(I)} C2^{(I)}(k,t) \right\} + \qquad (7)$$

$$+ \sum_{k=1}^{n} \left\{ b0_k^{(I)} S0^{(I)}(k,t) + b1_k^{(I)} S1^{(I)}(k,t) + b2_k^{(I)} S2^{(I)}(k,t) \right\}$$

Calculating the value of the spline $Se_{2,n}^{(I)}(t)$ and its derivatives at the nodes of the grid $\Delta_N^{(I)}$ and equating them with the corresponding value of trigonometric polynomials $T0_n^{(i)}(t)$, $T1_n^{(i)}(t)$, $T2_n^{(i)}(t)$, interpolating the function and its derivatives in the nodes of the same grid, it is easy to obtain systems of equations for the required coefficients $a0_0^{(I)}, a0_k^{(I)}, b0_k^{(I)}, a1_k^{(I)}, b1_k^{(I)}$ and $a2_k^{(I)}, b2_k^{(I)}$ of spline $Se_{2,n}^{(I)}(t)$. We confine ourselves to giving formulas for calculating the trigonometric Hermite spline that interpolates the function $f(t)$ and its derivatives up to and including the second order in the nodes of the grid $\Delta_N^{(0)}$; while still, we will omit the grid index. This spline has the form (7), where

$$C0(k,t) = \sum_{m=0}^{\infty} \frac{\cos(3mN+k)t}{(3mN+k)^4} \; ; \qquad S0(k,t) = \sum_{m=0}^{\infty} \frac{\sin(3mN+k)t}{(3mN+k)^4} \; ;$$

$$C1(k,t) = \sum_{m=0}^{\infty} \frac{\cos((3m+1)N+k)t}{((3m+1)N+k)^4} \; ; \qquad S1(k,t) = \sum_{m=0}^{\infty} \frac{\sin((3m+1)N+k)t}{((3m+1)N+k)^4} \; ;$$

$$C2(k,t) = \sum_{m=0}^{\infty} \frac{\cos((3m+2)N+k)t}{((3m+2)N+k)^4} \; ; \qquad S2(k,t) = \sum_{m=0}^{\infty} \frac{\sin((3m+2)N+k)t}{((3m+2)N+k)^4} \; ;$$

and coefficients $a0_0$, $a0_k$, $a1_k$, $a2_k$, $b0_k$, $b1_k$ $b2_k$ ($k = 1,2,...,n$) are found as solutions of system of equation

$$a0_0 = A0_0;$$

$$\begin{cases} a0_k \cdot 1 & + a1_k\, \alpha 10(k) + a2_k\, \alpha 20(k) = A0_k; \\ a0_k\, \alpha 01(k) + a1_k\, \alpha 11(k) + a2_k\, \alpha 21(k) = -B1_k; \\ a0_k\, \alpha 02(k) + a1_k\, \alpha 12(k) + a2_k\, \alpha 22(k) = -A2_k; \end{cases}$$

$$\begin{cases} b0_k \cdot 1 & + b1_k\, \alpha 10(k) + b2_k\, \alpha 20(k) = B0_k; \\ b0_k\, \alpha 01(k) + b1_k\, \alpha 11(k) + b2_k\, \alpha 21(k) = -A1_k; \\ b0_k\, \alpha 02(k) + b1_k\, \alpha 12(k) + b2_k\, \alpha 22(k) = -B2_k; \end{cases}$$

where $A0_0$, $A0_k$, $B0_k$, ($k = 1,2,...,n$)-coefficients of the trigonometric polynomial interpolating the function $f(t)$ in grid nodes $\Delta_N^{(0)}$, a $A1_k, B1_k$ and $A2_k, B2_k$ - coefficients of trigonometric polynomials interpolating derivatives, respectively 1 and 2 orders of magnitude $f(t)$ in nodes of the same grid;

$$g00(r,k) = \sum_{m=0}^{\infty} \frac{1}{(3mN+k)^{1+r}} \; ; \qquad g01(r,k) = \sum_{m=0}^{\infty} \frac{1}{(3mN+k)^r} \; ; \qquad g02(r,k) = \sum_{m=0}^{\infty} \frac{1}{(3mN+k)^{r-1}} \; ;$$

$$g10(r,k) = \sum_{m=0}^{\infty} \frac{1}{((3m+1)N+k)^{1+r}} \; ; \quad g11(r,k) = \sum_{m=0}^{\infty} \frac{1}{((3m+1)N+k)^r} \; ; \quad g12(r,k) = \sum_{m=0}^{\infty} \frac{1}{((3m+1)N+k)^{r-1}} \; ;$$

$$g20(r,k) = \sum_{m=0}^{\infty} \frac{1}{((3m+2)N+k)^{1+r}} \; ; \quad g21(r,k) = \sum_{m=0}^{\infty} \frac{1}{((3m+2)N+k)^r} \; ; \quad g22(r,k) = \sum_{m=0}^{\infty} \frac{1}{((3m+1)N+k)^{r-1}} \; ;$$

$$\alpha 00 = 1 \; ; \quad \alpha 01(r,k) = \frac{g01(r,k)}{g00(r,k)} \; ; \quad \alpha 02(r,k) = \frac{g02(r,k)}{g00(r,k)} \; ;$$

$$\alpha 10 = 1 \; ; \quad \alpha 11(r,k) = \frac{g11(r,k)}{g10(r,k)} \; ; \quad \alpha 12(r,k) = \frac{g12(r,k)}{g10(r,k)} \; ;$$

$$\alpha 20 = 1 \; ; \quad \alpha 21(r,k) = \frac{g21(r,k)}{g20(r,k)} \; ; \quad \alpha 22(r,k) = \frac{g22(r,k)}{g20(r,k)} \; .$$

Finally, it is advisable to make such a remark. When constructing Hermite trigonometric splines, we assumed that the values of the derivatives on the grid $\Delta_N^{(l)}$ are centered. However, this limitation is easy to remove, for example, in this way [12]. We consider a differential analogue of the Kronecker symbol

$$\frac{d^j}{dx^j} H_i = \begin{cases} 1, \ i = j; \\ 0, \ i \neq j, \end{cases} \quad i,j = 0,1,2.$$

Given this, the Hermit trigonometric spline, for example, for $p = 1$, we will look in the form

$$Se_{1,n}^{(0)}(t) = \frac{a_{0,0}^{(0)}}{2} H_0 + \frac{a_{1,0}^{(0)}}{2} H_1 +$$

$$+ \sum_{k=1}^{n} \left[ a0_k^{(0)} C0^{(0)}(k,t) + b0_k^{(0)} S0^{(0)}(k,t) + a1_k^{(0)} C1^{(0)}(k,t) + b1_k^{(0)} S1^{(0)}(k,t) \right],$$

where

$$a_{k,0}^{(l)} = \frac{1}{N} \sum_{j=1}^{N} f_j^{(l,k)} \; , \; (k = 0,1) \; .$$

In what follows, we will assume that when calculating the values of the Hermite trigonometric spline, we apply the operator to it. $\dfrac{d^0}{dx^0}$ ; with

$$\frac{d^0}{dx^0}H_0 = 1, \quad \frac{d^0}{dx^0}H_1 = 0.$$

and we have a formula (2).

When calculating the values of the derivative 1- we have order

$$\frac{d^1}{dx^1}H_0 = 0, \quad \frac{d^1}{dx^1}H_1 = 1,$$

and we obtain the formula for the derivative with a free term..

Thus, the proposed approach comes down to a solution 2-x  systems of algebraic equations of dimension $p+1$; in addition, known Fourier fast algorithms can be used to calculate the right-hand sides of the equations of these systems..

## CONCLUSIONS

1. A method for constructing Hermite trigonometric splines interpolating the function itself and its derivatives to $p$ -ro ( $p = 1, 2$ ) order inclusive in grid nodes $\Delta_N^{(l)}$, ( $N = 2n+1$, $n = 1, 2, ...$), for the implementation of which is necessary  solution 2-x systems of linear algebraic equations, each of which is a system $p+1$ order..

2. Thinking similarly, it is easy to construct trigonometric Hermite splines for larger values of the parameter $p$ .

3. The solution of the systems of linear algebraic equations necessary for the implementation of the proposed method can be done in advance; in practical application in each case, it is sufficient to substitute the coefficients of interpolation trigonometric polynomials into the solutions found..

4. When calculating the coefficients of interpolation trigonometric polynomials required for the construction of trigonometric Hermite polynomials by the proposed method, it is expedient to apply the known algorithms for fast Fourier transform.

5. Given the above, the proposed method of constructing trigonometric Hermite splines can be recommended for widespread use in practice.

## List of references